\documentclass[12pt]{article}
\usepackage{graphicx}
\usepackage{setspace}
\usepackage{bigints}
\usepackage{graphics} % for pdf, bitmapped graphics files
\usepackage{epstopdf}
\usepackage{epsfig} % for postscript graphics files
\usepackage{mathptmx} % assumes new font selection scheme installed
\usepackage{times} % assumes new font selection scheme installed
\usepackage{amsmath} % assumes amsmath package installed
\usepackage{amsthm}
\usepackage{amssymb}  % assumes amsmath package installed
\usepackage{algorithmicx}
\usepackage{algorithm}
\usepackage{bm}
\usepackage[noend]{algpseudocode}
\usepackage{cite}
\usepackage{subfig}
\usepackage{gensymb}
\usepackage{booktabs}
\usepackage{url}

\usepackage{tikz,pgfplots}
\pgfplotsset{compat=newest} 
\pgfplotsset{plot coordinates/math parser=false} 

\DeclareMathOperator*{\minimize}{minimize}

\DeclareMathOperator*{\subjectto}{subject\:to}

\DeclareMathAlphabet{\pazocal}{OMS}{zplm}{m}{n}

\newlength\figureheight 
\newlength\figurewidth  

\usepackage{algcompatible}
\usepackage{algpseudocode}

\makeatletter
\renewcommand*{\ALG@name}{Workflow}
\makeatother

\title{Optimization-based motion planning for multi-steered articulated vehicles}

\author
{\small Oskar Ljungqvist, Kristoffer Bergman and Daniel Axehill\\
	\\
	\small{Department of Automatic Control, Link\"oping University, Link\"oping, Sweden}\\
	\small{E-mail: \texttt{\{oskar.ljungqvist, kristoffer.bergman, daniel.axehill\}@liu.se}}
}

\date{}

\begin{document} 
		
\baselineskip 16pt
	
\maketitle 
%===============================================================================

\begin{abstract} % Abstract of not more than 250 words.
The task of maneuvering a multi-steered articulated vehicle in confined environments is difficult even for experienced drivers. In this work, we present an optimization-based trajectory planner targeting low-speed maneuvers in unstructured environments for multi-steered N-trailer vehicles, which are comprised of a car-like tractor and an arbitrary number of interconnected trailers with fixed or steerable wheels. The proposed trajectory planning framework is divided into two steps, where a lattice-based trajectory planner is used in a first step to compute a resolution optimal solution to a discretized version of the trajectory planning problem. The output from the lattice planner is then used in a second step to initialize an optimal control problem solver, which enables the framework to compute locally optimal trajectories that start at the vehicle's initial state and reaches the goal state exactly. The performance of the proposed optimization-based trajectory planner is evaluated in a set of practically relevant scenarios for a multi-steered 3-trailer vehicle with a car-like tractor where the last trailer is steerable.       
\end{abstract}

%%%%%%%%%%%%%%%%%%%%%%%%%%%%%%%%%%%%%%%%%%%%%%%%%%%%%%%%%%%%%%%
%% SECTION: INTRODUCTION
%
\section{Introduction}
In recent years, there has been a growing demand within the transportation sector to increase efficiency and reduce environmental impact related to transportation of both people and goods. This has lead to an increased interest for large capacity (multi-) articulated buses~\cite{michalek2019modular} and long tractor-trailer vehicle combinations~\cite{islam2015comparative}. In order to improve these long vehicles ability to maneuver in confined environments, some trailers (or wagons) can be equipped with steerable wheels. In the literature, these vehicles are commonly referred to as multi-steered N-trailer (MSNT) vehicles~\cite{orosco2002modeling}, which is a generalization of single-steered N-trailer (SSNT) vehicles. Compared to a SSNT vehicle, the additional steering capability enables an MSNT vehicle to be more flexible and agile, on the expense of increased difficulty in manually maneuvering the vehicle by a human driver. These difficulties partially arise due to the vehicle's increased degrees of freedom which are hard to successfully cope with for a human operator, and because of the specific kinematic and dynamic properties of an MSNT vehicle (see, e.g.,~\cite{tilbury1995multisteering,michalek2019modular,orosco2002modeling,islam2015comparative}).
To aid the driver, several advanced driver-assistance system concepts have been proposed to automatically steer the extra steerable wheel(s) in order to increase low-speed maneuverability or to reduce the so called off-tracking during tight cornering~\cite{odhams2011active,van2015active,michalek2019modular,varga2018robust}. 

Although a large amount of different motion planning techniques has been proposed in the literature for SSNT vehicles (see, e.g.,~\cite{sekhavat1997multi,hillary,evestedtLjungqvist2016planning,li2019trajectory,LjungqvistJFR2019}), there only exists a limited amount of work that consider the trajectory planning problem for special classes of MSNT vehicles (see e.g.~\cite{bushnell1995steering,tilbury1995multisteering,vidal2002real,Beyersdorfer2013tractortrailer,Yuan2017}). As a consequence, there is still a need to develop a trajectory planner that is able to solve the trajectory planning problem for a generic MSNT vehicle with car-like tractor that: i) can handle various state and input constraints, ii) allows a mixture of on-axle/off-axle hitched and steerable/non-steerable trailers, and iii) computes locally optimal trajectories by combining forward and backward motion.

The contribution of this work is a trajectory planning framework for an MSNT vehicle with car-like tractor targeting low-speed maneuvers in confined and unstructured environments. The framework extends some techniques presented in our previous work in~\cite{LjungqvistJFR2019} and is inspired by~\cite{bergman2019improved,bergman2019bimproved,bergman2018combining}. Here, a lattice-based trajectory planner is developed and used in a first step to compute a resolution optimal solution to a discretized version of the trajectory planning problem. The lattice planner uses a finite library of precomputed optimized maneuvers restricted to move the MSNT vehicle within a specified state-space discretization. In a second step, the output from the lattice planner is used to initialize a homotopy-based optimization step enabling the framework to compute locally optimal trajectories that starts at the vehicle's initial state and reaches desired goal states exactly. To the best of the authors knowledge, this paper presents the first trajectory planning framework for a generic MSNT vehicle with car-like tractor.

The reminder of the paper is organized as follows. In Section~\ref{p8:sec:model}, the kinematic vehicle model and the trajectory planning problem for the considered MSNT vehicle are presented, as well as an overview of the proposed trajectory planning framework. In Section~\ref{p8:sec:lattice_planner} and Section~\ref{p8:sec:homo_optimization}, the lattice-based trajectory planner and the homotopy-based optimization step are presented, respectively. Simulation results for an MS3T vehicle with car-like tractor are presented in Section~\ref{p8:sec:results} and the paper is concluded in Section~\ref{p8:sec:conclusions} by summarizing the contributions and discussing the directions for future work.

%%%%%%%%%%%%%%%%%%%%%%%%%%%%%%%%%%%%%%%%%%%%%%%%%%%%%%%%%%%%%%%
%% SECTION: Kinematic vehicle model and problem formulation
%
\section{Kinematic vehicle model and problem formulation}
\label{p8:sec:model}
\begin{figure}[t!]
	\centering
	\includegraphics[width=0.7\linewidth]{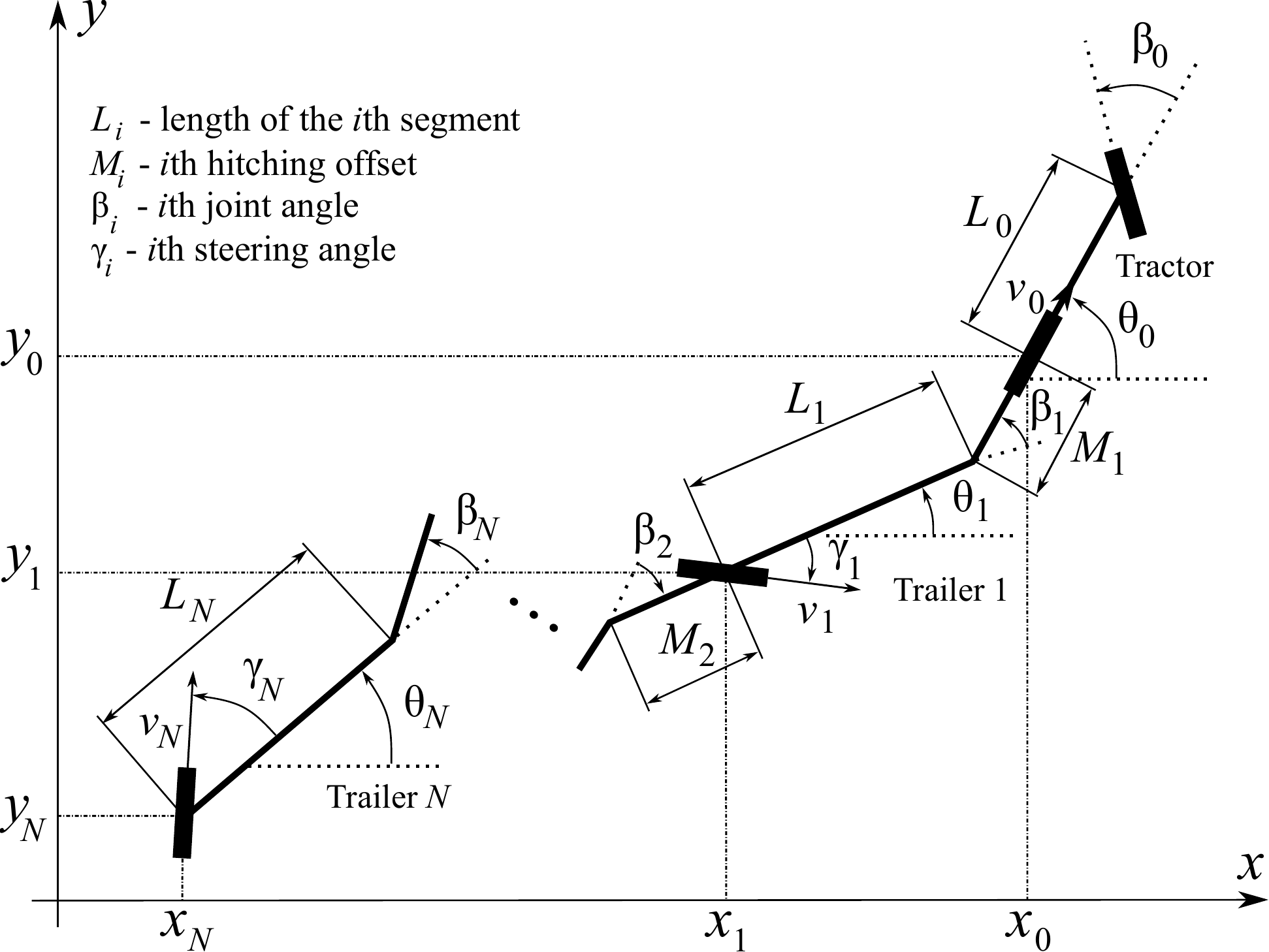}
	\caption{A schematic description of the geometric lengths and the vehicle configuration for a MSNT vehicle with car-like tractor in a global coordinate system (inspired and adapted from~\cite{michalek2019modular}).}
	\label{c8:fig:vehicle_model}%
\end{figure}

The MSNT vehicle with car-like tractor considered in this work is illustrated in Figure~\ref{c8:fig:vehicle_model}. The multi-body vehicle is composed of $N+1$ vehicle segments, including a car-like tractor and $N$ number of trailers that are equipped with steerable or non-steerable wheels. Each vehicle segment is characterized by a segment length $L_i>0$ and a signed hitching offset $M_i$. Since low-speed maneuvers are considered in this work, a kinematic vehicle model is used. The model is based on the work in~\cite{michalek2019modular} and is derived based on various assumptions such as the wheels are rolling without slipping and that the vehicle is operating on a flat surface. Moreover, it is assumed that the front wheel of the tractor is steerable and its rear wheel is non-steerable. The vehicle configuration consists of \mbox{$4 + N + S$} variables~\cite{michalek2019modular} where $S\in\{1,\hdots,N\}$ denotes the number of steerable trailers:
\begin{itemize}
	\item[--] the steering angle of the tractor's front wheel
	\begin{align}
		\beta_0 \in \pazocal Q_0=[-\bar\beta_0,\bar\beta_0], \quad \bar \beta_0 \in (0,\pi/2), 
				\label{c8:eq:front_steering_angle}
	\end{align}
	\item[--] the global position $(x_{N},y_{N})$ and orientation $\theta_{N}$ of the $N$th trailer in a fixed coordinate frame
	\begin{align}
		q_{N} = \begin{bmatrix}
			\theta_{N} & x_{N} & y_{N}
		\end{bmatrix}^T \in \mathbb S \times \mathbb R^2, \quad \mathbb S=(-\pi,\pi],
				\label{c8:eq:poes_trailer}
	\end{align} 
	\item[--] for $i=1,\hdots, N,$ a number of $N$ constrained joint angles
	\begin{align}
		\beta_i = \theta_{i-1} -\theta_{i} \in\pazocal B_i=[-\bar\beta_i,\bar\beta_i], \quad \bar \beta_i \in (0,\pi), 
		\label{c8:eq:joint_angles}
	\end{align}
	\item[--] and $S\in\{1,\hdots,N\}$ number of steering angles associated with steerable
	trailer wheels 
	\begin{align}
		\gamma_{s} \in \pazocal Q_{s}=[-\bar\gamma_{s},\bar\gamma_{s}], \quad \bar \gamma_{s} \in (0,\pi/2),
				\label{c8:eq:trailer_steering_angles}
	\end{align}
	where index $s\in \pazocal I_s \subseteq \{1,\hdots,N\}$ specifies which trailers that have steerable wheels. The configuration vector for the MSNT with car-like tractor will be defined as
	\begin{align}
		  q = \begin{bmatrix}
			\beta_0 & \beta_1 & \hdots & \beta_{N} & \bm \gamma_s^T  &  q_N^T 
		\end{bmatrix}^T \in \pazocal Q,
		\label{c8:eq:configuration_vector}
	\end{align}
	where $\bm \gamma_s$ represents a vector of trailer steering angles, and 
	$\pazocal Q =\pazocal Q_0\times \pazocal B_1 \hdots \times\pazocal B_N\times \underbrace{\pazocal Q_s\times \hdots\times \pazocal Q_s}_{S\text{-times}}\times \mathbb S\times \mathbb R^2$.
\end{itemize}
The leading car-like tractor is described by a kinematic single-track vehicle model and its orientation $\theta_0$ evolves as
\begin{align}
\dot \theta_0 = v_0\kappa_0(\beta_0),
\label{c8:eq:model_orientation_tractor}
\end{align}
where $\kappa_0(\beta_0) =\frac{\tan\beta_0}{L_0}$ is the curvature of the tractor and $v_0$ is the longitudinal velocity of its rear axle. The recursive formula for the transformation of the angular $\dot\theta_i$ and longitudinal $v_i$ velocities between any two neighboring vehicle segments are given by~\cite{michalek2019modular,orosco2002modeling}: 
\begin{align}
\begin{bmatrix}
\dot \theta_{i} \\  v_{i}
\end{bmatrix} &=
\underbrace{\begin{bmatrix}
-\cfrac{M_{i}}{L_{i}}\frac{\cos{(\beta_{i}-\gamma_{i})}}{\cos\gamma_{i}} & \frac{\sin{(\beta_{i}-\gamma_{i}+\gamma_{i-1})}}{L_{i}\cos\gamma_{i}} \\
M_{i}\frac{\sin{\beta_{i}}}{\cos\gamma_{i}} & \frac{\cos{(\beta_{i}+\gamma_{i-1})}}{\cos\gamma_{i}}
\end{bmatrix}}_{\triangleq J_i(\beta_{i},\gamma_i,\gamma_{i-1})}
\begin{bmatrix}
\dot \theta_{i-1} \\  v_{i-1}
\end{bmatrix}, \quad i=1,\hdots,N,
\label{c8:eq:velocity_transformation}
\end{align} 
where $\gamma_i$ denotes the steering angle of the $i$th trailer. For the car-like tractor, we have that $\gamma_0\equiv 0$ since its rear axle is non-steerable. Note that if the $j$th trailer is non-steerable, it suffices to take $\gamma_j=0$ in~\eqref{c8:eq:velocity_transformation}. 

Each trailer steering angle $\gamma_s$, $s\in\pazocal I_s$ and the tractors steering angle $\beta_0$ are modeled as double integrator systems
\begin{align}
\begin{split}
\dot \gamma_s &= \omega_s, \quad 
\dot \omega_s= u_{\omega_s}, \quad s\in\pazocal I_s,  \\
\dot \beta_0 &= \omega_0, \quad 
\dot \omega_0= u_{\omega_0},
\label{c8:eq:augumented_states}
\end{split}
\end{align}
where $\omega_0$, $\omega_s$, $s\in\pazocal I_s$ and $u_{\omega_0}$, $u_{\omega_s}$, $s\in\pazocal I_s$ denote steering angle rates and accelerations, respectively. This modeling is used to be able to penalize large rates and accelerations, and to enforce constraints in the form
\begin{align}
\begin{split}
\omega_s &\in \Omega_s =[-\bar \omega_s,\bar \omega_s], \quad u_{\omega_s} \in \pazocal U_s =[-\bar\Omega_s,\bar\Omega_s], \quad s\in\pazocal I_s, \\
\omega_0 &\in \Omega_0 =[-\bar \omega_0,\bar \omega_0], \quad u_{\omega_0} \in \pazocal U_0 =[-\bar\Omega_0,\bar\Omega_0],
\end{split}
\label{c8:eq:rate_acc_constraints_on_steering}
\end{align}
where the steering angle accelerations $u_{\omega_0}$ and \mbox{$u_{\omega_s}$, $s\in\pazocal I_s$} are treated as control signals. Similarly, the longitudinal velocity of the tractor $v_0$ is constrained as \mbox{$v_0\in \Omega_v = [-\bar v,\bar v]$} and its dynamics is modeled as a double integrator system 
\begin{align}
\label{c8:eq:long_dynamics}
\dot v_0 = a_0, \quad \dot a_0 = u_v
\end{align}
in order to respect constraints on the longitudinal acceleration $a_0\in \pazocal A = [-\bar a,\bar a]$ and jerk $u_v\in \pazocal U_v = [-\bar u_{v},\bar u_{v}]$. During the planning phase, the longitudinal jerk $u_{v}$ is treated as a control signal.

Since rolling without slipping of the wheels is assumed, the position of the last trailer evolves according to standard unicycle kinematics 
\begin{align}
\begin{split}
\dot x_{N} &= v_{N}\cos(\theta_{N}+\gamma_N), \\
\dot y_{N} &= v_{N}\cos(\theta_{N}+\gamma_N), 
\end{split}
\label{c8:eq:position_Ntrailer}
\end{align}
where its angular rate $\dot \theta_{N}$ and longitudinal velocity $v_{N}$ are given by
\begin{align}
\label{c8:eq:rate_speed_Ntrailer}
\begin{bmatrix}
\dot \theta_{N}\\
v_{N} 	
\end{bmatrix}=\prod_{i=0}^{N-1} \bm J_{N-i}(\beta_{N-i},\gamma_{N-i},\gamma_{N-i-1})\begin{bmatrix} v_0\kappa_0(\beta_0)\\
v_0 
\end{bmatrix},
\end{align}
which is derived by recursive usage of~\eqref{c8:eq:velocity_transformation} $N$ times together with~\eqref{c8:eq:model_orientation_tractor}. Combining~\eqref{c8:eq:position_Ntrailer} and \eqref{c8:eq:rate_speed_Ntrailer}, it is possible to compactly represent the model for the pose of the $N$th trailer $q_N  =[
\theta_{N} \hspace{5pt} x_{N} \hspace{5pt} y_{N}]^T$ as $\dot{ q}_N = v_0f_{ q_N}( q)$. In analogy, by introducing the vector $ c=[1\hspace{5pt} 0]^T$ together with~\eqref{c8:eq:model_orientation_tractor} and \eqref{c8:eq:velocity_transformation}, the time derivative of~\eqref{c8:eq:joint_angles} yields the joint-angle kinematics
\begin{align}
\dot\beta_i =& \dot \theta_{i-1} - \theta_i =  c^T \prod_{j=N-i-1}^{N-1}  J_{N-j}(\beta_{N-j},\gamma_{N-j},\gamma_{N-j-1})\begin{bmatrix}
v_0\kappa_0(\beta_0)\\
v_0 
\end{bmatrix} \nonumber \\
&- c^T\prod_{j=N-i}^{N-1}  J_{N-j}(\beta_{N-j},\gamma_{N-j},\gamma_{N-j-1})\begin{bmatrix} v_0\kappa_0(\beta_0)\\
v_0 
\end{bmatrix} \triangleq v_0f_{\beta_i}( q), \quad i=1,\hdots,N.
\label{c8:eq:joint_angle_kinematics}
\end{align}
Introduce the state vector as $ x = [
 q^T \hspace{5pt} \omega_0 \hspace{5pt} \bm\omega_s^T \hspace{5pt} v_0 \hspace{5pt} a_0
]^T\in\pazocal X$ and denote the control signal vector as $ u = [ u_{\omega_0}  \hspace{5pt}   \bm u_{\omega_s}^T  \hspace{5pt}  u_v
]^T\in\pazocal U$, where $\bm \omega_s$ and $\bm u_{\omega_s}$ represent vectors of trailer steering angle rates and accelerations, respectively, and 
\begin{align}
\begin{split}
\pazocal X  &= \pazocal Q\times\Omega_0\times\underbrace{\Omega_s\times\hdots\times\Omega_s}_{S\text{-times}}\times \Omega_v\times \pazocal A, \\
\pazocal U  &= \pazocal U_0\times \underbrace{\pazocal U_s\times\hdots\times\pazocal U_s}_{S\text{-times}}\times \pazocal U_v,
\end{split}
\label{c8:eq:pysicalConstraints}
\end{align} 
where $\text{dim}(\pazocal X) = 7+N+2S\triangleq n$ and \mbox{$\text{dim}(\pazocal U) = 2+S \triangleq m$}. The constraints in~\eqref{c8:eq:pysicalConstraints} will be referred to the vehicle's physical constraints arising from, e.g., actuator, mechanical or sensing limitations.
Finally, the kinematic model of the MSNT vehicle with car-like tractor is given in~\eqref{c8:eq:augumented_states},~\eqref{c8:eq:long_dynamics} and~\eqref{c8:eq:position_Ntrailer}--\eqref{c8:eq:joint_angle_kinematics}, which can compactly be represented as
\begin{align}
\dot{x}=f( x, u), \label{c8:vehicle_model}
\end{align}  
where $ f:\mathbb R^n\times \mathbb R^m\rightarrow \mathbb R^n$ is continuous and continuously differentiable with respect to $ x\in\pazocal X$ and $ u\in\pazocal U$. 
\subsection{Problem formulation}
The MSNT vehicle with car-like tractor is assumed to operate in a closed environment with only static obstacles $\pazocal X_{\text{obs}}$. The free-space where the vehicle is not in collision with any obstacle is defined as $\pazocal X_{\text{free}} = \pazocal X \setminus \pazocal X_{\text{obs}}$. Here, it is assumed that the obstacle-occupied region $\pazocal X_{\text{obs}}$ (hence also $\pazocal X_{\text{free}}$) can be described analytically, e.g., using different bounding regions~\cite{lavalle2006planning}. Since the free-space $\pazocal X_{\text{free}}$ is defined as the complement set of $\pazocal X_{\text{obs}}$, it is in general a non-convex set. 

The trajectory planning problem considered in this work is defined as follows: Compute a feasible and collision-free state and control signal trajectory $( x(t), u(t))$, $t\in[0,t_G]$ that moves the vehicle from its initial state $ x_I\in\pazocal X_{\text{free}}$ to a desired goal state $x_G\in\pazocal X_{\text{free}}$, while minimizing the cost functional $J$. This problem can be posed as a continuous-time optimal control problem (OCP) in the following general form 
\begin{align}
\begin{split}
\minimize_{u(\cdotp), \hspace{0.5ex}t_G }\hspace{3.7ex}
& J = t_G + \int_{0}^{t_G}l( x(t), u(t))\,\text{d}t 	\label{c8:eq:MotionPlanningOCP}\\
\subjectto\hspace{3ex}
& \dot{  x}(t) =  f( x(t), u(t)),  \\
& x(0) = x_I, \quad x(t_G) =  x_G,  \\ 
& x(t) \in \pazocal X_{\text{free}}, \quad
 u(t) \in \pazocal U,  
\end{split}
\end{align}
where the final time $t_G$ and $ u$ are optimization variables, and $l:\mathbb R^n\times \mathbb R^m\rightarrow \mathbb R_+$ is the cost function. In this work, the cost function is defined as
\begin{align}
l( x, u)=\lvert\lvert  x \rvert\rvert_{ Q}^2 + \lvert\lvert  u \rvert\rvert_{ R}^2,\label{c8:eq:costFunction}
\end{align}
where the weight matrices $Q\succeq 0$ and $R \succeq 0$. It is well-known that the OCP in~\eqref{c8:eq:MotionPlanningOCP} is in general hard to solve by directly invoking a state-of-the-art OCP solver~\cite{casadi,ipopt}. This is mainly because the vehicle model is nonlinear and the free-space $\pazocal X_{\text{free}}$ is in many applications non-convex. Hence, a proper initialization strategy for any OCP solver is often a necessity in order to converge to a good locally optimal (or even feasible) solution~\cite{bergman2019bimproved,zhang2018optimization}.

%%%%%%%%%%%%%%%%%%%%%%%%%%%%%%%%%%%%%%%%%%%%%%%%%%%%%%%%%%%%%%%%%%%%%%%%%%%%%%%%%%%
%% SECTION: Trajectory planner
%

\subsection{Trajectory planning framework}
\label{p8:sec:framework}
To efficiently and reliably solve the trajectory planning problem~\eqref{c8:eq:MotionPlanningOCP} for an MSNT vehicle with car-like tractor, we propose a framework that combines a lattice-based trajectory planner and an online optimization step. The framework is based on and extends the previous works in~\cite{LjungqvistJFR2019,bergman2019improved,bergman2019bimproved,bergman2018combining}. The extensions are made to account for the specific properties of an MSNT vehicle with car-like tractor. The general idea is that a lattice-based trajectory planner is used in a first step to compute an optimal solution to a discrete version of the trajectory planning problem~\eqref{c8:eq:MotionPlanningOCP} using a discretized state-space and a library of precomputed trajectories. The lattice planner is responsible for solving the combinatorial aspects of the trajectory planning problem, e.g., taking left or right around obstacles, and thus provides the latter optimization step with a proper initial guess. While keeping the combinatorial parts fixed, the objective of the optimization step is to further improve the initial guess computed by the lattice planner such that the resulting trajectory is a locally optimal solution to the original trajectory planning problem~\eqref{c8:eq:MotionPlanningOCP}. However, since the lattice planner uses a discretized state space, in general its computed trajectory does not satisfies the initial and goal state constraints in~\eqref{c8:eq:MotionPlanningOCP}. Therefore, the optimization step is also responsible for modifying the initial guess computed by the lattice planner such that the final optimized trajectory starts at the vehicle's initial state and reaches the goal state exactly. To handle this in a structured and numerically stable way, a homotopy-based optimization strategy is proposed that is inspired by the work in~\cite{bergman2018combining}. 

The main steps of the trajectory planning framework is summarized in Workflow~\ref{c8:alg1}, where the lattice planner \mbox{(Step 0 and 1)} and the optimization step (Step 2) are explained in detail in Section~\ref{p8:sec:lattice_planner} and Section~\ref{p8:sec:homo_optimization}, respectively.  
\vspace{5pt}
\begin{algorithm}[t!]
	\caption{The proposed trajectory planning framework for a MSNT vehicle with car-like tractor}
	\label{c8:alg1}
	\begin{description}
		\item \textbf{Step 0 (offline) -- State lattice construction:}
		\begin{enumerate}
			\item[a)] \textbf{\emph{State-space discretization:}} Specify the resolution of the discretized state-space $\pazocal X_d$
			\item[b)] \textbf{\emph{Motion primitive generation:}} Compute the set of motion primitives $\pazocal P$ by specifying a set of desired maneuvers and solve~\eqref{c8:eq:MotionPrimitiveGenOCP} using an OCP solver
			\item[c)] \textbf{\emph{Heuristic function:}} Precompute a HLUT by calculating the optimal cost-to-go on a grid in an obstacle-free environment 
		\end{enumerate}
		\item\textbf{Step 1 -- Online planning:}
		\begin{enumerate}
			\item[a)] \textbf{\emph{Initialization:}} Project the vehicle's initial state $ x_I$ and desired goal state $ x_G$ to $\pazocal X_d$
			\item[b)] \textbf{\emph{Graph search:}} Solve the discrete graph search problem~\eqref{c8:eq:OCP_discrete} using a graph search algorithm
		\end{enumerate}
		\item\textbf{Step 2 -- Homotopy-based optimization step:}
		\begin{enumerate}	
			\item[a)] \textbf{\emph{Initialization:}} Initialize the homotopy-based OCP solver with the solution computed by the lattice planner 
			\item[b)] \textbf{\emph{Optimization:}} Solve the relaxed trajectory planning problem~\eqref{c8:eq:MotionPlanningOCPhomo} using an OCP solver
			\item[c)] \textbf{\emph{Return:}} Send the computed solution to a trajectory-tracking controller or report failure
		\end{enumerate}
	\end{description}
\end{algorithm}
\vspace{10pt}

\section{Lattice-based trajectory planner}
\label{p8:sec:lattice_planner}
The idea with lattice-based trajectory planning is to restrict the solutions of the trajectory planning problem~\eqref{c8:eq:MotionPlanningOCP} to a lattice graph $\pazocal G = \langle \pazocal V, \pazocal E \rangle$, which is a graph embedded in an Euclidean space that forms a regular and repeated pattern~\cite{pivtoraiko2009differentially}. The lattice graph is constructed offline by discretizing the vehicle's state space $\pazocal X_d\subset\pazocal X$ and precompute a set of motion primitives $\pazocal P$.
Each vertex $v[k] \in \pazocal V$ is a vehicle state $ x[k]\in\pazocal X_d$ and each edge $ e_i\in \pazocal E$ represents of a motion primitive $ m_i\in\pazocal P$. A motion primitive is a feasible trajectory \mbox{$(x^i(t),u^i(t))$, $t\in[0,t_f^i]$} that moves the vehicle from an initial state $x[k]\in\pazocal X_d$ to a final state $x[k+1]\in\pazocal X_d$ while satisfying $x^i(\cdotp)\in\pazocal X$ and $u^i(\cdotp)\in\pazocal U$. A motion primitive is in this way designed to connect two vertices in the graph and the kinematic constraints~\eqref{c8:vehicle_model} and the physical constraints~\eqref{c8:eq:pysicalConstraints} are satisfied offline. The cardinality of the motion primitive set is $|\pazocal P|=M$ and the motion primitives that can be used from $x[k]$ is denoted $\pazocal P(x[k])\subseteq \pazocal{P}$. Moreover, since the MSNT vehicle is position invariant, the motion primitive set $\pazocal P$ can be computed from the position of the $N$th trailer at the origin. Each motion primitive $m_i$ can then be translated and reused for all other positions on the grid.

Let $x[k+1] = f_{p}(x[k],m_{i})$ denote the successor state when motion primitive $m_{i}\in\pazocal P$ is applied from $x[k]$ and denote $J_p(m_{i})$ as the stage-cost associated to this transition, which is given by
\begin{align}
J_p(m_i) = t^i_f + \int_{0}^{t^i_f}l(x^i(t),u^i(t))\text{d}t,
\end{align}
where $l(x^i,u^i)$ is defined in~\eqref{c8:eq:costFunction}.
The resulting trajectory taken by the vehicle when motion primitive \mbox{$m_i\in\pazocal{P}$} is applied from $x[k]$, is collision-free if it does not collide with any obstacle $c(m_i,x[k])\in\pazocal X_{\text{free}}$. 
Define \mbox{$u_p:\mathbb Z_+\rightarrow \{1,\hdots,M\}$} as a discrete and integer-valued decision variable that is selected by the lattice planner, where $u_p[k]$ specifies which motion primitive that is applied at stage $k$. Now, the continuous-time trajectory planning problem~\eqref{c8:eq:MotionPlanningOCP} is approximated by the following discrete-time OCP~\cite{LjungqvistJFR2019}:	
\begin{align}
\begin{split}
\minimize_{\{u_p[k]\}^{N_p-1}_{k=0}, \hspace{0.5ex} N_p}\hspace{3.7ex}
& J_{D} = \sum_{k=0}^{N_p-1}J_p( m_{u_p[k]}) \label{c8:eq:OCP_discrete} \\
\subjectto\hspace{3ex}
&  x[k+1] = f_{p}( x[k], m_{u_p[k]}),  \\
&  x[0] = \bar{x}_I, \quad  x[N_p] = \bar{ x}_G,  \\ 
&  m_{u_p[k]} \in \pazocal P( x[k]),  \\
& c( m_{u_p[k]}, x[k]) \in \pazocal X_{\text{free}}, \quad k=0,\hdots,N_p-1,
\end{split}
\end{align}
where $\bar{x}_I$ and $\bar{x}_G$ are obtained by projecting the actual initial state $ x_I$ and desired goal state $ x_G$ to the their closest neighboring state in $\pazocal X_d$. The decision variables to the problem in~\eqref{c8:eq:OCP_discrete} are the motion primitive sequence $\{u_p[k]\}^{N_p-1}_{k=0}$ and its length $N_p\in\mathbb Z_+$. 
A feasible solution is an ordered sequence of collision-free motion primitives $\{m_{u_p[k]}\}^{N_p-1}_{k=0}$, i.e., a trajectory $(x(t),u(t))$, $t\in [0, t_G]$ that connects the projected initial state $x(0)=\bar{x}_I$ with the projected goal state $ x(t_G)= \bar{x}_G$. Given the set of all feasible solutions to~\eqref{c8:eq:OCP_discrete}, an optimal solution is one that minimizes the cost function $J_{\text{D}}$. During online planning, the discrete-time OCP in~\eqref{c8:eq:OCP_discrete} can be solved using classical graph-search algorithms such as A$^*$ together with an informative precomputed free-space heuristic look-up table (HLUT) as heuristic function~\cite{knepper2006high}. A HLUT significantly reduces the online
planning time, as it takes the vehicle’s nonholonomic constraints into
account and enables perfect estimation of cost-to-go in free-space
scenarios with no obstacles.

\subsection{State-space discretization}
It is important that the resolution of the discretized state space $\pazocal X_d$ and the cardinality of the motion primitive set $\pazocal P$ are chosen such that the vehicle is sufficiently agile to maneuver in confined environments. However, as they also define the size of the lattice graph $\pazocal G$, both the resolution of $\pazocal X_d$ and the cardinality of $\pazocal P$ have to be chosen carefully in order to maintain a reasonable search time during online planning~\cite{pivtoraiko2009differentially}. Motivated by this, the position of the $N$th trailer $(x_N[k],y_N[k])$ is discretized to a uniform grid with resolution $r$ and its orientation is irregularly\footnote{$\Theta$ is the the set of unique angles $-\pi<\theta_{N}\leq \pi$ that can be generated by $\theta_{N} = \arctan2(i,j)$ for two integers $i,j\in\{-2,-1,0,1,2\}$.} discretized $\theta_N[k]\in\Theta$ into \mbox{$|\Theta|=16$} different orientations as proposed in~\cite{pivtoraiko2009differentially}. Additionally, the longitudinal velocity of the tractor is discretized as $v_0[k]\in V=\{-\bar v,0,\bar v\}$, where $\bar v$ is the tractor's maximal allowed speed. All remaining vehicle states are constrained to zero at each vertex in the graph, which implies that the MSNT vehicle is arranged in a straight configuration. This means that the joint angles $\beta_i[k]$, $i=1,\hdots, N$, the steering angles $\beta_0[k]$, $\bm\gamma_s[k]$, the steering angle rates $\omega_0[k]$, $\bm\omega_s[k]$ as well as the longitudinal acceleration of the tractor $a_0[k]$ are all constrained to zero at each $x[k]\in\pazocal X_d$. As a consequence, $\text{dim}(\pazocal X_d)=4$, since only the pose $ p_N[k]$ and the velocity of the tractor $v_0[k]$ are allowed to vary between different vertices in the graph. The proposed discretization will impose restrictions, but is motivated to enable fast online planning. 
%Conceptually, if more time is available, these constraints can be released. 
Moreover, since the output from the lattice planner will be used to warm start a second optimization step, it will improve the initial guess computed by the lattice planner such that the finally computed trajectory is a locally optimal solution to the original trajectory planning problem~\eqref{c8:eq:MotionPlanningOCP}. 

\subsection{Motion primitive generation}
The set of motion primitives $\pazocal P$ is computed offline by solving a finite
set of OCPs from all initial states $x_s^i\in\pazocal X_d$ with the position of the $N$th trailer at the origin, to a set of
final states $x_f^i\in\pazocal X_d$ in a bounded neighborhood in an obstacle-free
environment. This procedure can be performed manually as in~\cite{LjungqvistJFR2019} or using exhaustive search together with pruning strategies as proposed in~\cite{pivtoraiko2009differentially,CirilloIROS2014}. In both cases, the motion primitive generation procedure will become time consuming or requires a designer with high system knowledge. Therefore, here we use the maneuver-based motion primitive generation framework introduced in~\cite{bergman2019improved}. Instead of selecting pairs of discrete vehicle state to connect, a set of desired maneuvers from each initial state $x_s^i\in\pazocal X_d$ is selected and an OCP solver together with a rounding heuristic are used to automatically select the optimal final state $x_f^i\in\pazocal X_d$. Each maneuver is defined with a terminal manifold in the form $g^i(x^i(t_f^i)) = 0$ where $g:\mathbb R^n\rightarrow \mathbb R^l$ and $l<n$, where $n-l$ is the degrees of freedom for the terminal state constraint. To compute a maneuver-specified motion primitive $m_i\in\pazocal P$, the following continuous-time OCP is first solved  
\begin{align}
\begin{split}
\minimize_{u^i(\cdotp), \hspace{0.5ex}t^i_p }\hspace{3.7ex}
& J_p(m_i) 	\label{c8:eq:MotionPrimitiveGenOCP}\\
\subjectto\hspace{3ex}
& \dot{x}^i(t) = f(x^i(t),u^i(t)),  \\
&x^i(0) =x^i_s, \quad g^i(x^i(t_f^i)) = 0,  \\ 
&x^i(t) \in \pazocal X^i, \quad
u^i(t) \in \pazocal U^i.  
\end{split}
\end{align}
Here it is not required that $x^i(t_f^i)\in\pazocal X_d$. To ensure that the final state $x_f^i=x^i(t_f^i)\in\pazocal X_d$, a rounding heuristic is used and the closest
neighboring states represented in the discretized state-space $\pazocal X_d$
from $x^i(t_f^i)$ are evaluated and the solution with lowest objective
functional value is selected as the resulting motion primitive $m_i$. Finally, since the MSNT with car-like tractor is orientation invariant, rotational symmetries are exploited to reduce the number of OCPs needed to be solved~\cite{pivtoraiko2009differentially}. For more details of the motion primitive generation framework, the reader is referred to~\cite{bergman2019improved}. Note that the vehicle's physical constraints $\pazocal X^i\subseteq \pazocal X $ and $\pazocal U^i \subseteq \pazocal U $ in~\eqref{c8:eq:MotionPrimitiveGenOCP} are defined to be maneuver dependent, which is not the case in~\cite{bergman2019improved}. This extension is made to automatically generate similar maneuvers, i.e., same terminal manifold $g^i(x^i(t_f^i))=0$, but resulting in different optimal final states $x^i_f$ and final times $t_f^i$. This additional freedom can be used to design a more flexible lattice planner or, e.g., to adapt to a change in available trailer steering angles $\gamma_s$, $s\in\pazocal I_s$ during different maneuvers.

As proposed in~\cite{bergman2019improved}, the motion primitive set is build upon optimized straight, heading change and parallel maneuvers. The heading change and parallel maneuvers are only possible to use from states where the tractor has nonzero velocity, i.e. $v_{0,s}=\pm \bar v$, and are designed to end up in the same final velocity $v_{0,f}=v_{0,s}$. Additionally, short straight maneuvers from $v_{0,s}\in V$ to some $v_{0,f}\in V$ are also optimized to enable the tractor to reduce, increase and keep a constant longitudinal velocity. The heading change and parallel maneuvers are computed using the following terminal manifolds.

\begin{figure}[t!]
 	\centering
 	\includegraphics[width=0.55\linewidth]{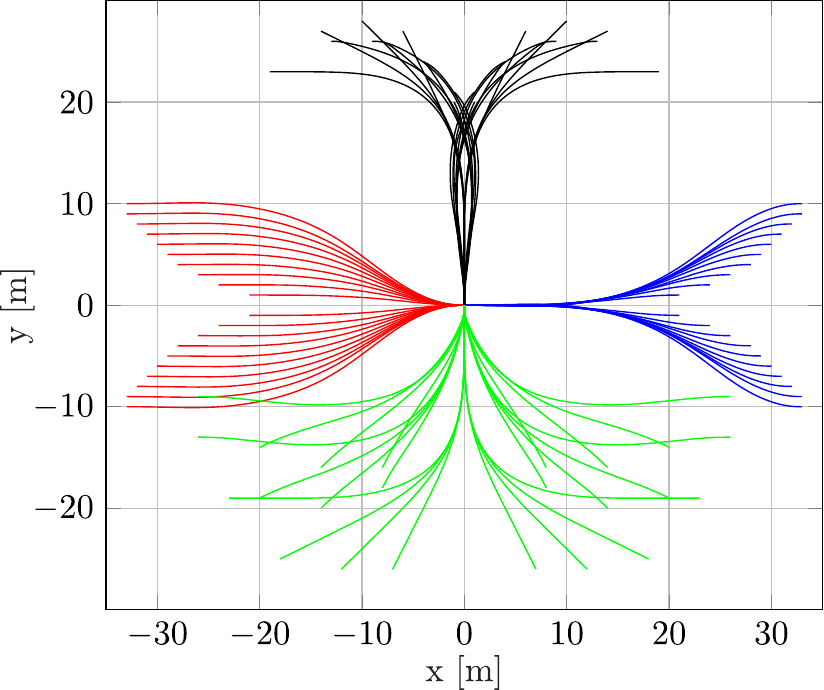}
 	%\begin{tikzpicture}
 	%\node[anchor=south west] (myplot) at (0,0) {
 	%	\input{iFAC_controlMulti/Figures/primitives.tex}
 	%};
 	%\end{tikzpicture}
 	\caption{A subset of the motion primitives in $\pazocal P$ for an MS3T vehicle with a car-like tractor from initial position at the origin $(x_{3,s},y_{3,s})=(0,0)$ to different final states $x_f\in\pazocal X_d$. The colored lines are the trajectories in $(x_3(\cdotp),y_3(\cdotp))$ for the different maneuvers. The set of heading change maneuvers from $(\theta_{3,s},v_{0,s})=(\pi/2,1)$ (black) and from $(\theta_{3,s},v_{0,s})=(\pi/2,-1)$ (green). The set of parallel maneuvers from $(\theta_{3,s},v_{0,s})=(0,1)$ (blue) and from $(\theta_{3,s},v_{0,s})=(0,-1)$ (red).}
 	\label{c8:fig:primitives}
 \end{figure}

\vspace{5pt}
\textbf{Heading change maneuvers:} By specifying the vehicle's physical constraints $\pazocal X^i\subseteq \pazocal X$ and $\pazocal U^i \subseteq \pazocal U $, a heading change maneuver from an initial state $x^i_s\in\pazocal X_d$ with pose \mbox{$p^i_{N,s}=[
0\hspace{5pt} 0\hspace{5pt} \theta^i_{N,s}]^T$} and $v^i_{0,s} = \pm \bar v$, to a user-defined adjacent orientation $\theta^i_{N,f}\in\Theta \setminus \theta^i_{N,s}$ is optimized by solving~\eqref{c8:eq:MotionPrimitiveGenOCP} using the following terminal constraint
\begingroup
\renewcommand*{\arraystretch}{1.3}
\begin{align}
g^i\left(p^i_N(t^i_f),v^i_{0,f}\right) = \begin{bmatrix}
\theta^i_N(t^i_f) - \theta^i_{N,f} \\
v^i_{0}(t^i_f) - v^i_{0,s}
\end{bmatrix}=0,
\label{c8:eq:heading_change_constraint}
\end{align}
\endgroup
which implies that $x^i_N(t_f)$ and $y^i_N(t_f)$ are free variables for the OCP solver to select. Note that the vehicle states that are left out from the argument to $g^i$ are all constrained to zero to guarantee that $x^i_f\in\pazocal X_d$. Examples of computed heading change maneuvers from $(\theta_{3,s},v_{0,s})=(\pi/2,\pm 1)$ are depicted in Figure~\ref{c8:fig:primitives} for an MS3T vehicle with car-like tractor where the last trailer has steerable wheels, i.e., $\pazocal I_s = \{3\}$. Here, the allowed trailer steering angle $|\gamma_3|\leq \bar\gamma_3$ is alternated using $\bar \gamma_3 = 0$, $0.175$ and $0.35$ rad, resulting in different types of optimal trajectories.
 %\tikzexternaldisable
 %\tikzexternalenable
 
\vspace{5pt}
\textbf{Parallel maneuvers:} A parallel maneuver from initial state $x^i_s\in\pazocal X_d$ with pose $p^i_{N,s}=[
0\hspace{5pt} 0\hspace{5pt} \theta^i_{N,s}]^T$ and \mbox{$v^i_{0,s} = \pm \bar v$}, is defined with a user-defined lateral displacement $z^i_{\text{lat}}$ in $(x^i_{N,f},y^i_{N,f})$ with respect to the initial orientation $\theta^i_{N,s}$. This maneuver can be optimized by solving~\eqref{c8:eq:MotionPrimitiveGenOCP} using the following terminal constraint 
%\begingroup\makeatletter\def\f@size{7.6}\check@mathfonts
\begin{align}
g^i\left(p^i_N(t^i_f),v^i_{0,f}\right) = \begin{bmatrix}
y^i_N(t^i_f)\cos\theta^i_{N,s} + x^i_N(t^i_f)\sin\theta^i_{N,s} - z^i_{\text{lat}}\\
\theta^i_N(t^i_f) - \theta^i_{N,s} \\
v^i_{0}(t^i_f) - v^i_{0,s}
\end{bmatrix}=0.
\label{c8:eq:parallel_constraint}
\end{align}
%\endgroup 
Here, the final position of the $N$th trailer $(x_N(t_f^i),y_N(t_f^i))$ is restricted to a line defined by the first row in $\eqref{c8:eq:parallel_constraint}$.      
Examples of computed parallel maneuvers for an MS3T vehicle with car-like tractor using \mbox{$\bar\gamma_3=0.35$ rad} from \mbox{$(\theta_{3,s},v_{0,s})=(0,\pm 1)$} can be seen in Figure~\ref{c8:fig:primitives}.

When the motion primitive set $\pazocal P$ is computed, a free-space heuristic look-up table (HLUT) is computed using techniques presented in~\cite{knepper2006high}. The HLUT is computed offline by solving several obstacle-free graph-search problems~\eqref{c8:eq:OCP_discrete} from all initial states $\bar {x}_I\in\pazocal X_d$ with $(x_{N,I},y_{N,I})=(0,0)$, to final states \mbox{$\bar {x}_G\in\pazocal X_d$} with positions on a bounded grid around the origin. This computation can be done efficiently using Dijkstra's search, as the optimal cost-to-come is simply recorded and stored in the HLUT~\cite{knepper2006high}. Moreover, in analogy to the motion primitive
generation, the size of the HLUT is kept small by exploiting the position and orientation invariance properties of $\pazocal P$~\cite{CirilloIROS2014}. 

%During online planning, the HLUT is used as heuristic function to guide the used graph-search algorithm. As shown in~\cite{knepper2006high}, a HLUT significantly reduces the online
%planning time, as it takes the vehicle’s nonholonomic constraints into
%account and enables perfect estimation of cost-to-go in free-space
%scenarios with no obstacles.

\section{Homotopy-based optimization step}
\label{p8:sec:homo_optimization}
Similar to~\cite{bergman2019bimproved}, the optimization step is used to improve the initial guess $(x(t),  u(t))$, \mbox{$t\in[0, t_G]$} computed by the lattice planner such that the final trajectory $( x^*(t),  u^*(t))$, \mbox{$t\in[0,t^*_G]$} is a locally optimal solution to~\eqref{c8:eq:MotionPlanningOCP}. Since the lattice planner uses a discretized state space $\pazocal X_d$, in general its computed state trajectory does not satisfy the initial and goal state constraints in~\eqref{c8:eq:OCP_discrete}. Thus, the optimization step should not only improve the initial guess but also make it feasible in the original problem formulation~\eqref{c8:eq:MotionPlanningOCP}. To handle this in a structured way, a homotopy-based initialization strategy is used that is inspired by the work in~\cite{bergman2018combining}. The idea is to start from a relaxed problem that is easy to solve and then gradually transform the relaxed problem to the original one. Here, these ideas are applied on the initial and goal state constraints in~\eqref{c8:eq:MotionPlanningOCP} such that the solution obtained from the lattice planner is feasible to the relaxed problem. By letting $\epsilon_p^T = [\epsilon_{p,I}\hspace{5pt}\epsilon_{p,G}]\in[0,1]^2$ denote the homotopy parameters~\cite{bergman2018combining}, the initial and goal state constraints in~\eqref{c8:eq:MotionPlanningOCP} are relaxed to 
\begin{align}\label{c8:eq:homoInitGoal}
\begin{split}
	& x(0) = \epsilon_{p,I}{\bar x}_I + (1-\epsilon_{p,I}){ x}_I, \\
	& x(t_G) = \epsilon_{p,G}{\bar x}_G + (1-\epsilon_{p,G}){ x}_G. 
\end{split}
\end{align}  
When $\epsilon_p^T = [1\hspace{5pt} 1]$, the initial guess from the lattice planner is feasible to the relaxed version of~\eqref{c8:eq:MotionPlanningOCP} and when $\epsilon_p^T = [0\hspace{5pt} 0]$, the original problem in~\eqref{c8:eq:MotionPlanningOCP} is obtained. One possibility is to start with $\epsilon_p^0 = [1\hspace{5pt} 1]^T$ and repeatedly solve the relaxed version of~\eqref{c8:eq:MotionPlanningOCP} using an OCP solver where $\epsilon_p^k$ is gradually decreased using a fixed step-size $\Delta\epsilon_p$ until $ \epsilon_p^{k} = [0\hspace{5pt} 0]^T$ is reached~\cite{bergman2018combining}. In this work, the idea is instead to let the OCP solver automatically modify the homotopy parameters using a penalty method~\cite{nocedal2006numerical}. Define the linear penalty as $ c_p^T\epsilon_p$, where $ c_p\in\mathbb R^2_{++}$ and define the relaxed version of the trajectory-planning problem~\eqref{c8:eq:MotionPlanningOCP} as    
\begin{align}
\begin{split}
\minimize_{ u(\cdotp), \hspace{0.5ex}t_G,\hspace{0.5ex}\epsilon_p }\hspace{3.7ex}
& J_{H} = J +  c_p^T\epsilon_p 	\label{c8:eq:MotionPlanningOCPhomo}\\
\subjectto\hspace{3ex}
& \dot{ x}(t) =  f( x(t), u(t)),  \\
& x(0) = \epsilon_{p,I}{\bar x}_I + (1-\epsilon_{p,I}){ x}_I, \\
& x(t_G) = \epsilon_{p,G}{\bar x}_G + (1-\epsilon_{p,G}){ x}_G, \\
& x(t) \in \pazocal X_{\text{free}}, \quad
 u(t) \in \pazocal U, \quad  \epsilon_p\in[0,1]^2,
\end{split}
\end{align}
which is initialized with the solution from the lattice planner and $\epsilon_p^T = [1\hspace{5pt} 1]$. By choosing $c_p$ sufficiently large, the OCP solver will automatically adjust the step size of $\epsilon_p$ and converge to \mbox{$\epsilon_p^T = [0\hspace{5pt} 0]$} if a feasible solution to~\eqref{c8:eq:MotionPlanningOCP} exists in the homotopy class selected by the lattice planner~\eqref{c8:eq:OCP_discrete}~\cite{bergman2018combining}. As previously mentioned, if \mbox{$\epsilon_p^* = [0\hspace{5pt} 0]^T$} is obtained, a locally optimal solution $( x^*(t), u^*(t))$, $t\in[0,t^*_G]$ to the original trajectory planning problem~\eqref{c8:eq:MotionPlanningOCP} is obtained which can then be sent to a trajectory-tracking controller for plan execution.  

Note that if $c_p$ is not chosen sufficiently large, a solution with \mbox{$\epsilon_p^T = [0\hspace{5pt} 0]$} may not be obtained even though one exists~\cite{nocedal2006numerical}. In that case, one may need to increase $c_p$ and continue the solution process. However, in extensive simulation trials presented in Section~\ref{p8:sec:results}, it is shown that the proposed homotopy-based optimization step is able to reliably compute locally optimal solutions to~\eqref{c8:eq:MotionPlanningOCPhomo} with \mbox{$\epsilon_p^T = [0\hspace{5pt} 0]$} without modifying $c_p$ in all problem instances. 

\section{Simulation results}
\label{p8:sec:results}
In this section, the proposed trajectory planning framework is evaluated in two complicated parking problem scenarios for an MS3T with car-like tractor where only trailer $N=3$ is steerable, i.e., $\pazocal I_s = \{3\}$, and a mixture of off-axle ($M_1\neq 0$) and on-axle hitching ($M_2=M_3=0$). Using the recursive model presented in Section~\ref{p8:sec:model} it is now straightforward to derive the vehicle model~\eqref{c8:vehicle_model} for this specific MS3T vehicle with configuration $ q=[\beta_0\hspace{5pt}\beta_1\hspace{5pt}\beta_2\hspace{5pt}\beta_3\hspace{5pt}\gamma_3\hspace{5pt}\theta_3\hspace{5pt}x_3\hspace{5pt}y_3]^T$, augmented state vector $ x = [ q^T\hspace{5pt}\omega_0\hspace{5pt}\omega_3\hspace{5pt}v_0\hspace{5pt}a_0]^T$ and control signals $ u = [u_{\omega_0}\hspace{5pt}u_{\omega_3}\hspace{5pt}u_{v}]^T$. The vehicle model $\dot{ x}= f( x, u)$ is presented in Appendix A and the values of the vehicle's parameters used in this section are summarized in Table~\ref{c8:tab:vehicle_parameters}.
\begin{table}[t!]
	\caption{Vehicle parameters for the MS3T vehicle.}
	\centering
	\begin{tabular}{l l}
		\hline \noalign{\smallskip} Vehicle parameter  & Value   \\  \hline \noalign{\smallskip}	
		Tractor's wheelbase $L_0$        &   4.6 m  \\ 
		Length of the off-hitch $M_1$        &   1.6 m  \\
		Length of trailer 1 $L_1$            &   2.5 m  \\ 
		Length of trailer 2 $L_2$            &   7.0 m  \\
		Length of trailer 3 $L_3$            &   7.0 m  \\
		Maximum joint angles $\bar\beta_i$, $i=1,2,3$ & 0.87 rad \\
		Maximum steering angle tractor $\bar\beta_0$ & 0.73 rad \\   
		Maximum steering-angle rate tractor $\bar\omega_0$ & $0.8$ rad/s\\   
		Maximum steering-angle acceleration tractor $\bar\Omega_0$ & $10$ rad/s$^2$\\  
		Maximum steering angle trailer 3  $\bar\gamma_3$ & 0.35 rad \\   
		Maximum steering-angle rate trailer 3 $\bar\omega_3$ & $0.4$ rad/s\\   
		Maximum steering-angle acceleration trailer 3 $\bar\Omega_3$ & $10$ rad/s$^2$\\  
		Maximum longitudinal speed tractor $\bar v$ &1 m/s \\
		Maximum longitudinal acceleration tractor $\bar a$ &$1$ m/s$^2$ \\
		Maximum longitudinal jerk tractor $\bar u_{v}$ &40 m/s$^2$ \\
		\hline \noalign{\smallskip}
	\end{tabular}
	\label{c8:tab:vehicle_parameters}
\end{table}
The cost function is chosen as
\begin{align}
l( x, u) = \frac{1}{2}\left(\beta_0^2+\gamma_s^2+10\omega_0^2+10\omega_3^2 + a_0^2 +  u^T u\right),
\end{align} 
which is used in all steps of the trajectory planning framework. The linear cost on the homotopy parameters in the optimization step is chosen as $c_p^T = [1000\hspace{5pt}1000]$. The lattice planner is implemented in \texttt{C++}, whereas the motion primitive generation and the homotopy-based optimization step are both implemented in Python using \textrm{CasADi}~\cite{casadi}, where \textrm{IPOPT} is used as nonlinear programming problem solver. All simulations are performed
on a laptop computer with an Intel Core i7-4600U@2.1GHz CPU.

The motion primitive set consists of heading change, parallel and straight trajectories where and a subset of the motion primitive set \mbox{$\pazocal P_{\text{MS3T}}$ ($|\pazocal P_{\text{MS3T}}|=2080$)} can be seen in Figure~\ref{c8:fig:primitives}. From each initial heading with nonzero longitudinal velocity, there are $20$ parallel and $24$ heading change maneuvers. The heading change maneuvers are computed using three different limits on the trailer steering angle $\bar \gamma_3=0$, $0.175$ and $0.35$ rad, respectively, and the parallel maneuvers using \mbox{$\bar \gamma_3=0.35$ rad}. After the motion primitive set is computed, a free-space HLUT is computed on a square grid $80\times80$ m centered around the origin. To evaluate if the trajectory planner is able to exploit the additional trailer steering, it is compared with an SS3T vehicle, i.e., $\bar \gamma_3=0$, with the same vehicle parameters and the difference that only $8$ heading change maneuvers exist in the motion primitive set $\pazocal P_{\text{SS3T}}$ ($|\pazocal P_{\text{SS3T}}|=1124$). 

%
%
%\tikzexternaldisable
\begin{figure}[t!]
	\centering
	\includegraphics[width=0.8\linewidth]{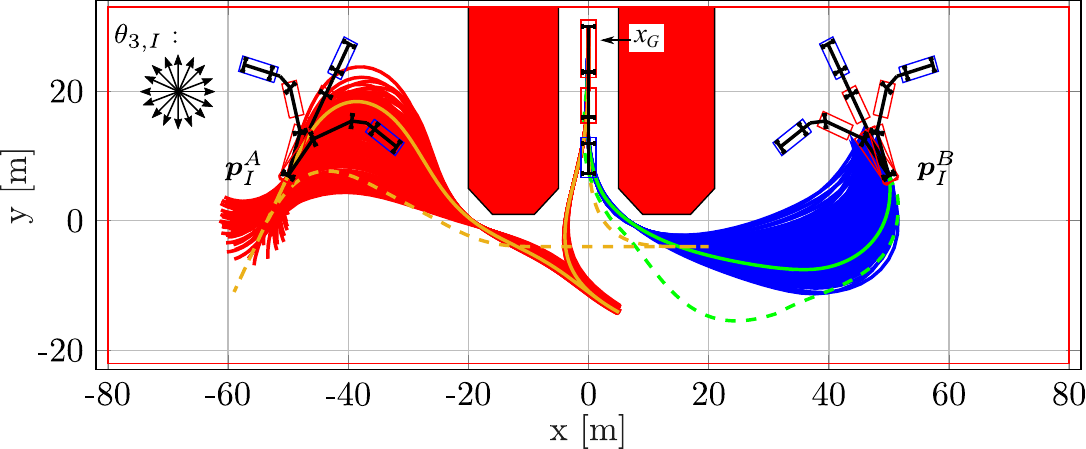}
	%	\begin{tikzpicture}
	%	\node[anchor=south west] (myplot) at (0,0) {
	%		\input{Figures/loading_bay_ex.tex}
	%	};
	%\begin{scope}[x={(myplot.south east)}, y={(myplot.north west)}]
	%\node[text=black] at (0.25,0.64) {\small $\bm p_I^A$};
	%\node[text=black] at (0.85,0.64) {\small $\bm p_I^B$};
	%\node[text=black] at (0.17,0.9) {\small $\theta_{3,I}:$};
	%\end{scope}
	%\end{tikzpicture}
	\caption{Loading-bay parking scenario using the proposed trajectory planner for SS3T and MS3T from two initial positions ($ p_I^A$ and $ p_I^B$), symmetric with respect to the loading bay, initial orientations $\theta_{3,I}\in\Theta$ and zero initial joint angles, to the goal state $ x_G$. The solutions for the position ($x^*_3(\cdotp),y^*_3(\cdotp)$) from 1029 perturbed initial states for SS3T ($ p_I^A$) and MS3T ($p_I^B$) are displayed by blue (red) solid lines for MS3T (SS3T). The initial guess computed by the lattice planner for ($ x_3(\cdotp), y_3(\cdotp)$) and optimization step for ($x^*_3(\cdotp),y^*_3(\cdotp)$) from straight configuration is marked with green (orange) dashed line and green (orange) solid line, respectively, for MS3T (SS3T).}
	\label{c8:fig:Loading_bay}
\end{figure}
%\tikzexternalenable

\begin{table}[b!]
	\caption{Results from loading-bay parking scenario in Figure~\ref{c8:fig:Loading_bay} for 32 problems. $\bar t_{\text{lat}}$ and $\bar t_{\text{ocp}}$  is average computation time for lattice planner and optimization step, respectively. $\bar J_D$ and $\bar J_H$ is average objective functional value for the solutions from lattice planner and optimization step, respectively. $\bar r_{\text{imp}}$ and $\bar t_{G,\text{imp}}$ is average cost and time improvement between the lattice planner's and the optimization step's solutions, respectively.}
	\centering
	\begin{tabular}{l l l l l l l}
		\hline \noalign{\smallskip} Vehicle & $\bar t_{\text{lat}}$[s] & $\bar t_{\text{ocp}}$[s] & $\bar J_D$ & $\bar J_H$ & $\bar r_{\text{imp}}$ & $\bar t_{G,\text{imp}}$[s]\\  \hline \noalign{\smallskip}	
		SS3T  & 0.11 & 9.3 &  170.8   & 128.4   &   -25\%   & -26.8  \\ 
		MS3T  & 0.09 & 3.0 &  126.8   & 100.8   &   -21\%  & -14.4 \\
		\hline \noalign{\smallskip}
	\end{tabular}
	\label{c8:tab:loading_bay_results}
\end{table}

The first planning scenario is a loading-bay parking problem that is illustrated in Figure~\ref{c8:fig:Loading_bay}. The obstacles and vehicle bodies are described by bounding circles~\cite{lavalle2006planning}, where in total, the vehicle bodies are described using 8 bounding circles of radius 2 m. The objective of the trajectory planner is to plan a trajectory from 32 different initial states $x_I\in\pazocal X_d$ (see Figure~\ref{c8:fig:Loading_bay}) to the goal state $x_G$. One solution example is provided for the lattice planner (dashed line) and optimization step (solid line) for SS3T (orange) and MS3T (green), respectively, for symmetric planning problems. The results show that the planned trajectory for MS3T is purely in backward motion, as apposed to SS3T which needs to combine forward and backward motion due to less steering capability. A summary of the simulation results are provided in Table~\ref{c8:tab:loading_bay_results}. The average computation time for the lattice planner is only {$0.1$ s} for both SS3T and MS3T, whereas the optimization step takes in average {$3.0$ s} for MS3T and {$9.0$ s} for SS3T. However, a signification reduction of both average cost and time improvement of the solutions are obtained. Note that the average time improvement of the solutions computed by the optimization step $\bar t_{G,\text{imp}}$ is significantly larger compared to the optimization step's average computation time $\bar t_{G,\text{ocp}}$. Thus, when the optimization step is added, the combined average computation and execution time becomes significantly lower.

%\tikzexternaldisable
\begin{figure}[t!]
	\centering
	\includegraphics[width=1\linewidth]{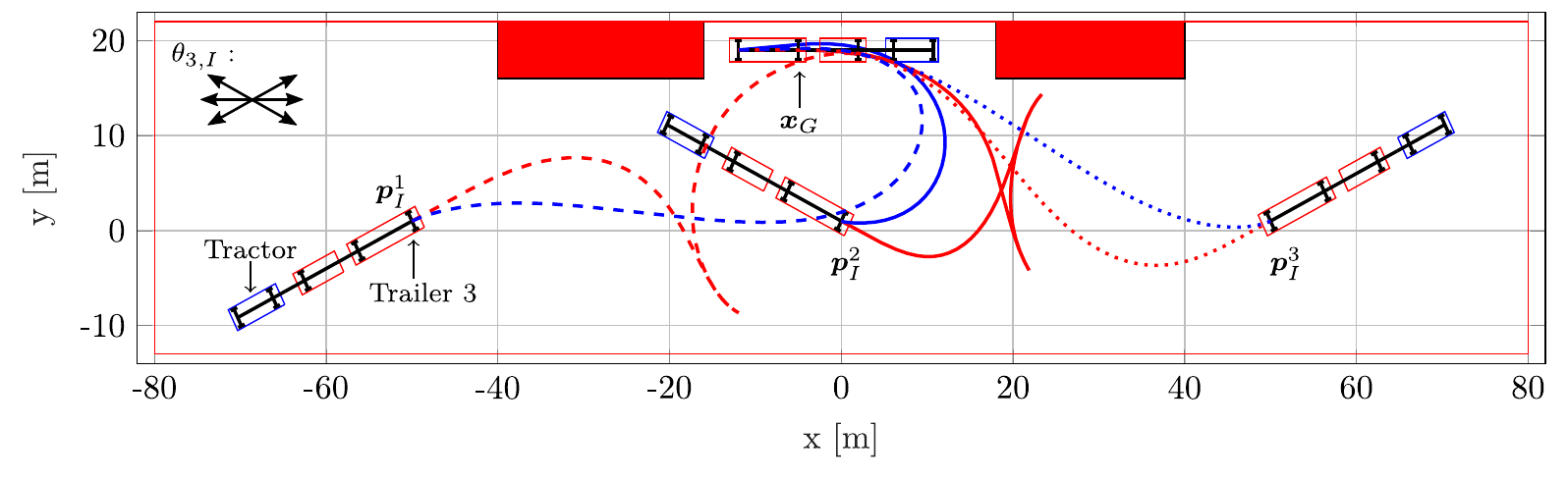}
	%\begin{tikzpicture}
	%\node[anchor=south west] (myplot) at (0,0) {
	%	\input{Figures/parallel_parking_ex.tex}
	%};
	%\begin{scope}[x={(myplot.south east)}, y={(myplot.north west)}]
	%\node[text=black] at (0.51,0.74) {\small $\bm x_G$};
	%\draw[->] (0.51,0.78) to [out=270,in=270] (0.51,0.85);  
	%\node[text=black] at (0.16,0.48) {\small Tractor};
	%\draw[->] (0.16,0.45) to [out=90,in=90] (0.16,0.39);  
	%\draw[->] (0.264,0.425) to [out=270,in=270] (0.264,0.50); 
	%\node[text=black] at (0.27,0.39) {\small Trailer 3};
	%\node[text=black] at (0.25,0.60) {\small $\bm p_I^1$};
	%\node[text=black] at (0.54,0.45) {\small $\bm p_I^2$};
	%\node[text=black] at (0.82,0.45) {\small $\bm p_I^3$};
	%\node[text=black] at (0.13,0.86) {\small $\theta_{3,I}:$};
	%\end{scope}
	%\end{tikzpicture}
	\caption{Parallel parking scenario using the proposed trajectory planner for SS3T and MS3T from three different initial positions ($ p^1_I$, $ p^2_I$ and $ p^3_I$), six different initial orientations $\theta_{3,I}$ and zero initial joint angles, to the goal state $ x_G$. The blue (red) lines illustrate the planned trajectories after the optimization step for the position of trailer 3 $(x^*_3(\cdotp),y^*_3(\cdotp))$ for MS3T (SS3T) from three selected initial states.}
	\label{c8:fig:parallel_parking}
\end{figure}
%\tikzexternalenable

\begin{table}[b!]
	\caption{Results from parallel parking scenario in Figure~\ref{c8:fig:parallel_parking} for 18 problems. See Table~\ref{c8:tab:loading_bay_results} for description of the variables.}
	\centering
	\begin{tabular}{l l l l l l l}
		\hline \noalign{\smallskip} Vehicle & $\bar t_{\text{lat}}$[s] & $\bar t_{\text{ocp}}$[s] & $\bar J_D$ & $\bar J_H$ & $\bar r_{\text{imp}}$ & $\bar t_{G,\text{imp}}$[s]\\  \hline \noalign{\smallskip}	
		SS3T  & 5.6 & 8.4 &  270.6   & 122.4   &   -55\%   & -109.6  \\ 
		MS3T  & 11.3 & 2.4 &  207.9   & 98.1   &   -52\%  & -80.4 \\
		\hline \noalign{\smallskip}
	\end{tabular}
	\label{c8:tab:parallel_parking_results}
\end{table}

To evaluate how the homotopy-based optimization step handles an infeasible initialization, 1029 perturbed initial states for both SS3T and MS3T are evaluated (see Figure~\ref{c8:fig:Loading_bay}), where blue and red trajectories are related to MS3T and SS3T, respectively. The initial joint angles are perturbed with \mbox{$\beta_i=-30^\circ,-20^\circ,\hdots,30^\circ$}, $i=1,2,3$ and initial orientation with $\theta_{3,I}=-10^\circ,0^\circ,10^\circ$. In all cases, the used optimization step is able to handle the infeasible initial guess obtained from the lattice planner, i.e., the value of homotopy parameter $\epsilon_{p,I}^*=0$ in all cases. That is, in all cases, a solution to the original trajectory planning problem~\eqref{c8:eq:MotionPlanningOCP} is obtained. Moreover, the average computation time of the optimization step is \mbox{$3.4$ s} for MS3T and \mbox{$26.2$ s} for SS3T. Hence, the active trailer steering also reduces the computation load of the OCP solver.  

\begin{figure}[t!]
	\vspace{-10pt}
	\centering
	\setlength\figureheight{0.22\columnwidth}
	\setlength\figurewidth{0.85 \columnwidth} 
	\captionsetup[subfloat]{captionskip=-3pt} 
	\subfloat[][The joint angle between tractor and trailer 1 $\beta_1(\cdotp)$ (black), joint angle between the trailer 1 and trailer 2 $\beta_2(\cdotp)$ (blue), and joint angle between the trailer 2 and trailer 3 $\beta_3(\cdotp)$ (red). Their limits are displayed by dashed black lines.]{
		\includegraphics[width=0.8\linewidth]{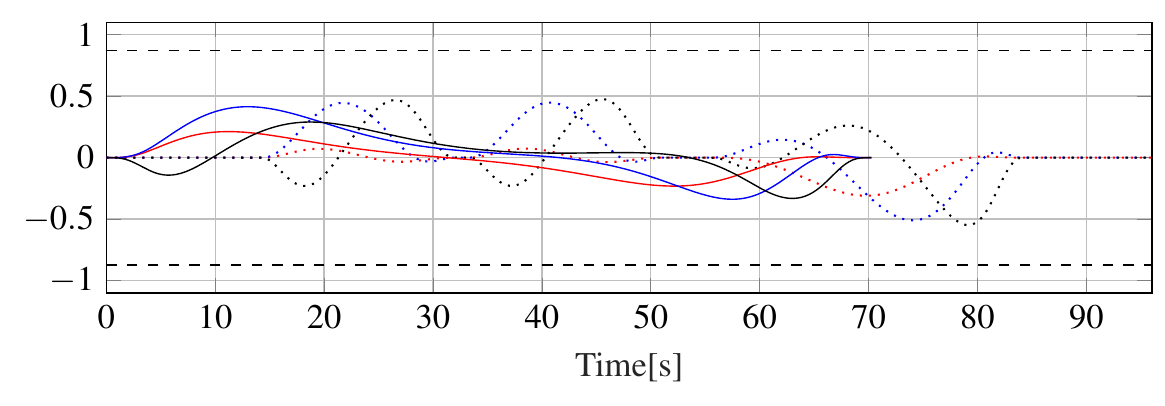}
		%	\begin{tikzpicture}
		%\node[anchor=south west] (myplot) at (0,0) {
		%	\input{Figures/parallel_parking_ex_joint_angles.tex}
		%};
		%\end{tikzpicture}
		\label{p8:fig:parallel_parking_joint_angles}
	}
	\quad
	\subfloat[][The steering angle of the tractor $\beta_0$ (blue), velocity $v_0$ (black) and steering angle of trailer 3 $\gamma_3$ (red). Their limits are displayed by dashed-dotted lines.]{
		\includegraphics[width=0.8\linewidth]{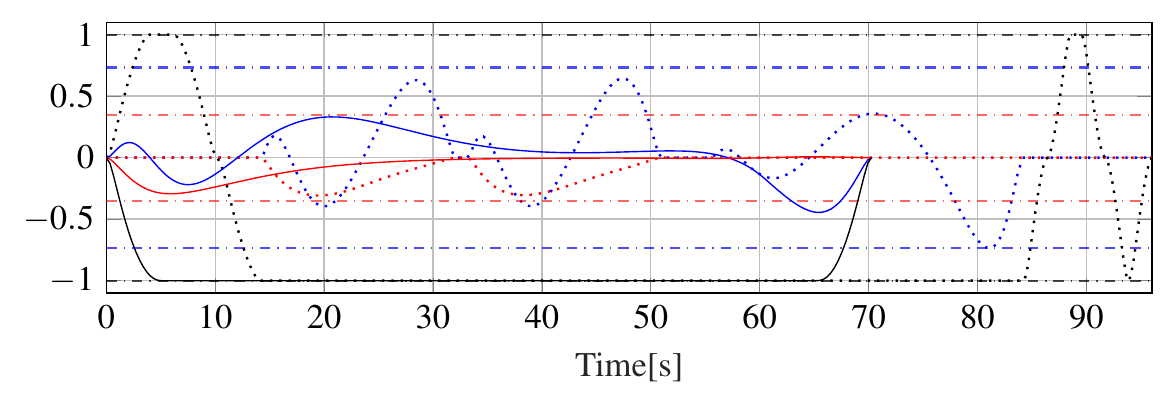}
		%\begin{tikzpicture}
		%\node[anchor=south west] (myplot) at (0,0) {
		%\input{Figures/parallel_parking_ex_controls.tex}
		%};
		%\end{tikzpicture}
		\label{p8:fig:parallel_parking_controls}
	}	
	\caption{A subset of the resulting trajectories for the parallel parking scenario in Figure~\ref{c8:fig:parallel_parking} from initial state $p_I^3$ with MS3T for optimized (solid lines) and lattice initial guess (dotted lines).}
	\label{p8:fig:parallel_parking_state_traj}
\end{figure}

The second planning scenario is a parallel parking problem with 18 different problems that is illustrated in Figure~\ref{c8:fig:parallel_parking} and the results are summarized in Table~\ref{c8:tab:parallel_parking_results}. This scenario is a confined environment which affects the average computation time of the lattice planner $\bar t_{\text{lat}}$, which is \mbox{$11.3$ s} for MS3T and \mbox{$5.6$ s} for SS3T. This is because the HLUT is drastically underestimating the cost-to-go in this confined environment. Therefore, both the average cost improvement $\bar r_{\text{imp}}$ (MS3T: $-52\%$, SS3T: $-55\%$) and time improvement $\bar t_{G,\text{imp}}$ (MS3T: $-80.4$~s, SS3T: $-109.6$~s) of the optimization step are significant. The confined environment does however not affect the average computation time of the optimization step which is 2.4 s for MS3T and \mbox{8.4 s} for SS3T. Moreover, as can be seen in the three highlighted planning problems in Figure~\ref{c8:fig:parallel_parking}, the final optimized solutions for the MS3T only needs to reverse, as opposed to the SS3T which needs to combine forward and backward motion in two cases. Finally, Figure~\ref{p8:fig:parallel_parking_state_traj} shows the difference between the trajectories from the lattice planner and the optimization step for the selected planning problem from position $p_I^3$ in Figure~\ref{c8:fig:parallel_parking}. As can be seen, the trajectories for the two steering angles, the longitudinal velocity and the joint angles are significantly smoother after the optimization step, at the same time as the final time is decreased from $96$ s to $70$ s. 

%\tikzexternaldisable
%\tikzexternalenable

\newpage
\section{Conclusions}
\label{p8:sec:conclusions}
An optimization-based trajectory planner for multi-steered articulated vehicles is proposed that targets low-speed maneuvers in unstructured environments. The proposed trajectory planner is divided into two steps, where a lattice planner is used in a first step to compute an optimal solution to a discretized version of the trajectory planning problem using a library of precomputed trajectories. 
In a second step, the output from the lattice planner is then used to initialize a homotopy-based optimization step, which enables the framework to compute a locally optimal solution that starts at the vehicle's initial state and reaches the goal state exactly. The performance of the proposed optimization-based trajectory planner is evaluated in a set of practically relevant scenarios for a multi-steered 3-trailer vehicle where the last trailer is steerable. In the simulations, it is shown that the framework can solve challenging trajectory planning problems and that the proposed optimization step provides a significant improvement in terms of reduced objective functional value and final time, at the same time as it enables the framework to plan from and to a larger set of different vehicle states. 

As future work we would like to develop a trajectory-tracking controller and evaluate the framework in real-world experiments on a full-scale test vehicle.       

\section*{Appendix A}
In this section, the model~\eqref{c8:vehicle_model} for the specific MS3T vehicle used in the simulation trails is presented. The matrices $ J_1$, $ J_2$ and $ J_3$ describing the longitudinal and angular velocity transformations between neighboring vehicle segments~\eqref{c8:eq:velocity_transformation} are 
\begin{align}
	\label{c8:MS3T_vel_matrices}
	\begin{split}
		 J_1(\beta_1,0,0)&=\begin{bmatrix}
			-\frac{M_{1}}{L_{1}}\cos\beta_{1} & \frac{\sin\beta_1}{L_{1}} \\[10pt]
			M_{1}\sin\beta_1 & \cos\beta_{1}
		\end{bmatrix}, 
	\quad 
		 J_2(\beta_2,0,0)=\begin{bmatrix}
			0 & \frac{\sin\beta_2}{L_2} \\[10pt]
			0 & \cos\beta_2
		\end{bmatrix}, \\
	 J_3(\beta_3,\gamma_3,0)&=\begin{bmatrix}
		0 & \frac{\sin(\beta_3-\gamma_3)}{L_3\cos\gamma_3} \\[10pt]
		0 & \frac{\cos\beta_3}{\cos\gamma_3}
	\end{bmatrix},
	\end{split}
\end{align}
since $M_2=M_3=0$ and $\gamma_0=\gamma_1=\gamma_2=0$. Thus, by inserting~\eqref{c8:MS3T_vel_matrices} in~\eqref{c8:eq:rate_speed_Ntrailer} the longitudinal velocity $v_{3}$ of trailer 3 can be written as
\begin{align}
	v_{3} = v_0\frac{\cos\beta_3}{\cos\gamma_3}\cos\beta_2\left(\frac{M_1}{L_0}\sin\beta_1\tan\beta_0 + \cos\beta_1\right).
	\label{c8:MS3T_vel_trailer3}
\end{align} 
Finally, using~\eqref{c8:MS3T_vel_matrices} and~\eqref{c8:MS3T_vel_trailer3} in~\eqref{c8:eq:position_Ntrailer}--\eqref{c8:eq:joint_angle_kinematics}, the complete model~\eqref{c8:vehicle_model} for the MS3T vehicle becomes
\begin{align}
	\begin{split}
		\label{c8:vehicle_model_MS3T}
	\dot\beta_1 =& v_0\left(\frac{\tan\beta_0}{L_0} -\frac{\sin\beta_1}{L_1} + \frac{M_1}{L_0L_1}\cos\beta_1\tan\beta_0\right),\\
	\dot\beta_2 =& v_0\left[\left(\frac{\sin\beta_1}{L_1} - \frac{M_1}{L_0L_1}\cos\beta_1\tan\beta_0\right) \right. \\ &\left.
	-\frac{\sin{\beta_{2}}}{L_{2}}\left(\frac{M_1}{L_0}\sin\beta_1\tan\beta_0 + \cos\beta_1\right)\right], \\
	\dot\beta_3 =& v_0\left[\frac{\sin{\beta_{2}}}{L_{2}}\left(\frac{M_1}{L_0}\sin\beta_1\tan\beta_0 + \cos\beta_1\right) \right. \\ &\left. - \frac{\sin{(\beta_{3}-\gamma_{3})}}{L_{3}\cos\gamma_{3}}\cos\beta_2\left(\frac{M_1}{L_0}\sin\beta_1\tan\beta_0 + \cos\beta_1\right)\right], \\
	\dot \theta_{3} =& 
	v_0\frac{\sin{(\beta_{3}-\gamma_{3})}}{L_{3}\cos\gamma_{3}}\cos\beta_2\left(\frac{M_1}{L_0}\sin\beta_1\tan\beta_0 + \cos\beta_1\right), \\
	\dot x_3 =&v_0\frac{\cos\beta_3}{\cos\gamma_3}\cos\beta_2\left(\frac{M_1}{L_0}\sin\beta_1\tan\beta_0 + \cos\beta_1\right) \cos(\theta_{3}+\gamma_3), \\
	\dot y_3 =&  v_0\frac{\cos\beta_3}{\cos\gamma_3}\cos\beta_2\left(\frac{M_1}{L_0}\sin\beta_1\tan\beta_0 + \cos\beta_1\right) \sin(\theta_{3}+\gamma_3), \\
	\dot \beta_0 =& \omega_0, \quad \dot \omega_0 = u_{\omega_0}, \\
	\dot \gamma_3 =& \omega_3,\quad	\dot \omega_3 = u_{\omega_3}, \\
	\dot v_0 =& a_0, \quad \dot a_0 = u_v.
	\end{split}
\end{align}

\bibliography{root}

\begin{thebibliography}{10}

\bibitem{casadi}
J.~A.~E. Andersson et~al.
\newblock {CasADi} -- {A} software framework for nonlinear optimization and
  optimal control.
\newblock {\em Mathematical Programming Computation}, 2018.

\bibitem{bergman2018combining}
K.~Bergman and D.~Axehill.
\newblock Combining homotopy methods and numerical optimal control to solve
  motion planning problems.
\newblock In {\em Proceedings of the 2018 IEEE Intelligent Vehicles Symposium},
  pages 347--354, 2018.

\bibitem{bergman2019improved}
K.~Bergman, O.~Ljungqvist, and D.~Axehill.
\newblock Improved optimization of motion primitives for motion planning in
  state lattices.
\newblock In {\em Proceedings of the 2019 IEEE Intelligent Vehicles Symposium},
  2019.

\bibitem{bergman2019bimproved}
K.~Bergman, O.~Ljungqvist, and D.~Axehill.
\newblock Improved path planning by tightly combining lattice-based path
  planning and optimal control.
\newblock \emph{Under review for possible publication in IEEE Transactions on
  intelligent vehicles}. Pre-print available at {arXiv}:
  \url{https://arxiv.org/abs/1903.07900}, 2019.

\bibitem{Beyersdorfer2013tractortrailer}
S.~Beyersdorfer and S.~Wagner.
\newblock Novel model based path planning for multi-axle steered heavy load
  vehicles.
\newblock In {\em Proceedings of the 16th International Conference on
  Intelligent Transportation Systems}, pages 424--429, Oct 2013.

\bibitem{bushnell1995steering}
L.~G. Bushnell et~al.
\newblock Steering three-input nonholonomic systems: the fire truck example.
\newblock {\em The International Journal of Robotics Research}, 14(4):366--381,
  1995.

\bibitem{CirilloIROS2014}
M.~Cirillo, T.~Uras, and S.~Koenig.
\newblock A lattice-based approach to multi-robot motion planning for
  non-holonomic vehicles.
\newblock In {\em Proceedings of the 2014 IEEE/RSJ International Conference on
  Intelligent Robots and Systems}, pages 232--239, 2014.

\bibitem{evestedtLjungqvist2016planning}
N.~Evestedt, O.~Ljungqvist, and D.~Axehill.
\newblock Motion planning for a reversing general 2-trailer configuration using
  {C}losed-{L}oop {RRT}.
\newblock In {\em Proceedings of the 2016 IEEE/RSJ International Conference on
  Intelligent Robots and Systems}, pages 3690--3697, 2016.

\bibitem{islam2015comparative}
M.~M. Islam et~al.
\newblock A comparative study of multi-trailer articulated heavy-vehicle
  models.
\newblock {\em Proceedings of the Institution of Mechanical Engineers, Part D:
  Journal of Automobile Engineering}, 229(9):1200--1228, 2015.

\bibitem{knepper2006high}
R.~A. Knepper and A.~Kelly.
\newblock High performance state lattice planning using heuristic look-up
  tables.
\newblock In {\em Proceedings of the 2006 IEEE/RSJ International conference on
  Intelligent Robots and Systems}, pages 3375--3380, 2006.

\bibitem{hillary}
F.~Lamiraux et~al.
\newblock Motion planning and control for hilare pulling a trailer.
\newblock {\em IEEE Transactions on Robotics and Automation}, 15(4):640--652,
  Aug 1999.

\bibitem{lavalle2006planning}
S.~M. LaValle.
\newblock {\em Planning algorithms}.
\newblock Cambridge {U}niversity {P}ress, 2006.

\bibitem{li2019trajectory}
B.~Li et~al.
\newblock Trajectory planning for a tractor with multiple trailers in extremely
  narrow environments: A unified approach.
\newblock In {\em Proceeding of the 2019 International Conference on Robotics
  and Automation}, pages 8557--8562, 2019.

\bibitem{LjungqvistJFR2019}
O.~Ljungqvist et~al.
\newblock A path planning and path-following control framework for a general
  2-trailer with a car-like tractor.
\newblock {\em Journal of Field Robotics}, 36(8):1345--1377, 2019.

\bibitem{michalek2019modular}
M.~M. Michalek.
\newblock Modular approach to compact low-speed kinematic modelling of
  multi-articulated urban buses for motion algorithmization purposes.
\newblock In {\em Proceeding of the 2019 IEEE Intelligent Vehicles Symposium},
  pages 2060--2065, 2019.

\bibitem{nocedal2006numerical}
J.~Nocedal and S.~Wright.
\newblock {\em Numerical optimization}.
\newblock Springer Science \& Business Media, 2006.

\bibitem{odhams2011active}
A.~Odhams et~al.
\newblock Active steering of a tractor--semi-trailer.
\newblock {\em Proceedings of the Institution of Mechanical Engineers, Part D:
  Journal of Automobile Engineering}, 225(7):847--869, 2011.

\bibitem{orosco2002modeling}
R.~Orosco-Guerrero et~al.
\newblock Modeling and dynamic feedback linearization of a multi-steered
  n-trailer.
\newblock {\em IFAC Proceedings Volumes}, 35(1):103--108, 2002.

\bibitem{pivtoraiko2009differentially}
M.~Pivtoraiko et~al.
\newblock Differentially constrained mobile robot motion planning in state
  lattices.
\newblock {\em Journal of Field Robotics}, 26(3):308--333, 2009.

\bibitem{sekhavat1997multi}
S.~Sekhavat et~al.
\newblock Multilevel path planning for nonholonomic robots using semiholonomic
  subsystems.
\newblock {\em The International Journal of Robotics Research}, 17(8):840--857,
  1998.

\bibitem{tilbury1995multisteering}
D.~Tilbury et~al.
\newblock A multisteering trailer system: conversion into chained form using
  dynamic feedback.
\newblock {\em IEEE Transactions on robotics and automation}, 11(6):807--818,
  1995.

\bibitem{van2015active}
N.~Van De~Wouw et~al.
\newblock Active trailer steering for robotic tractor-trailer combinations.
\newblock In {\em Proceeding of the 54th IEEE Conference on Decision and
  Control}, pages 4073--4079, 2015.

\bibitem{varga2018robust}
B.~Varga et~al.
\newblock Robust tracking controller design for active dolly steering.
\newblock {\em Proceedings of the Institution of Mechanical Engineers, Part D:
  Journal of Automobile Engineering}, 232(5):695--706, 2018.

\bibitem{vidal2002real}
T.~Vidal-Calleja, M.~Velasco-Villa, and E.~Aranda-Bricaire.
\newblock Real-time obstacle avoidance for trailer-like systems.
\newblock In {\em Proceeding of the 3th International Symposium on Robotics and
  Automation}, 2002.

\bibitem{ipopt}
A.~W\"achter and L.~T. Biegler.
\newblock On the implementation of a primal-dual interior point filter line
  search algorithm for large-scale nonlinear programming.
\newblock {\em Mathematical Programming}, 106(1):25--57, 2006.

\bibitem{Yuan2017}
J.~{Yuan}.
\newblock Hierarchical motion planning for multisteering tractor-trailer mobile
  robots with on-axle hitching.
\newblock {\em IEEE/ASME Transactions on Mechatronics}, 22(4):1652--1662, Aug
  2017.

\bibitem{zhang2018optimization}
X.~Zhang et~al.
\newblock Autonomous parking using optimization-based collision avoidance.
\newblock In {\em Proceedings of the 2018 IEEE Conference on Decision and
  Control}, pages 4327--4332, Dec 2018.

\end{thebibliography}
\bibliographystyle{abbrv}

\end{document}